\documentclass{amsart}
\usepackage[dvips]{epsfig}
\usepackage{amsfonts}
\usepackage{amsmath}
\usepackage{amssymb}
\usepackage{graphics}
\newtheorem{thm}{Theorem}[subsection]

 \theoremstyle{definition}

 \numberwithin{equation}{subsection}

\begin{document}
\title[Positive Dehn Twist Expressions for Some Elements of Finite Order]
 {Positive Dehn Twist Expressions for Some Elements of Finite Order in the Mapping Class Group}

\author{ Yusuf Z. Gurtas }

\address{Department of Mathematics, Suffolk CCC, Selden, NY, USA}

\email{gurtasy@sunysuffolk.edu}

\thanks{}




\dedicatory{}

\commby{}

\begin{abstract}
Positive Dehn twist products for some elements of finite order in
the mapping class group of a 2-dimensional closed, compact,
oriented surface $\Sigma_g$, which are rotations of $\Sigma_g$
through $2\pi /p$, are presented. The homeomorphism invariants of
the resulting simply connected symplectic 4- manifolds are
computed.

\end{abstract}

\maketitle

\section*{Introduction}
Positive Dehn twist expressions for a new set of involutions in
the mapping class group of 2-dimensional closed, compact, oriented
surfaces were presented in \cite{Gu1}. In a later article
 the idea was extended  to bounded surfaces and from
that, the  positive Dehn twist expressions for involutions on
closed surfaces that are obtained by joining several copies of
bounded surfaces together were also obtained, \cite{Gu2}. It is
the purpose of this article to present positive Dehn twist
expressions for some elements of finite order in $M_g$, which are
rotations through $2\pi /p$ on the surface $\Sigma_g$. This is
shown in section 2 using the expressions for the involutions
described in \cite{Gu1} . In section 3 the homeomorphism
invariants of the resulting simply connected symplectic 4-
manifolds that are realized as Lefschetz fibrations are computed.
\section{Review}

\begin{figure}[htbp]
     \centering  \leavevmode
     \psfig{file=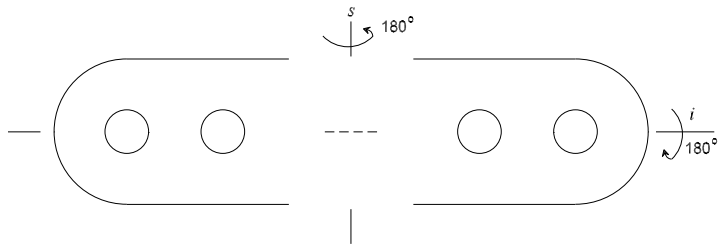,width=3.750in,clip=}
     \caption{The vertical and horizontal involutions}
     \label{twoinvolutions.fig}
 \end{figure}

Let $i$ represent the hyperelliptic (horizontal) involution and
$s$ represent the vertical involution as shown in Figure
\ref{twoinvolutions.fig}.

\begin{figure}[htbp]
     \centering  \leavevmode
     \psfig{file=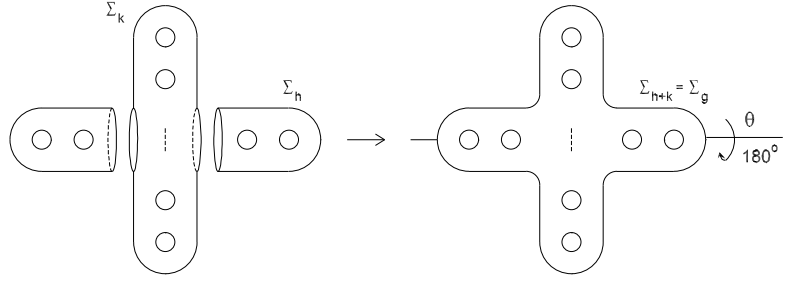,width=4.50in,clip=}
     \caption{The involution $\theta$ on the surface $\Sigma_{h+k}$}
     \label{gluingtwoinvolutions.fig}
 \end{figure}

If $i$ is the horizontal involution on a surface $\Sigma_{h}$ and
$s$ is the vertical involution on a surface $\Sigma_{k}$,
$k-$even, then let $\theta$ be the horizontal involution on the
surface $\Sigma_{g}$, where $g=h+k$, obtained as in Figure
\ref{gluingtwoinvolutions.fig}.

\begin{figure}[htbp]
     \centering  \leavevmode
     \psfig{file=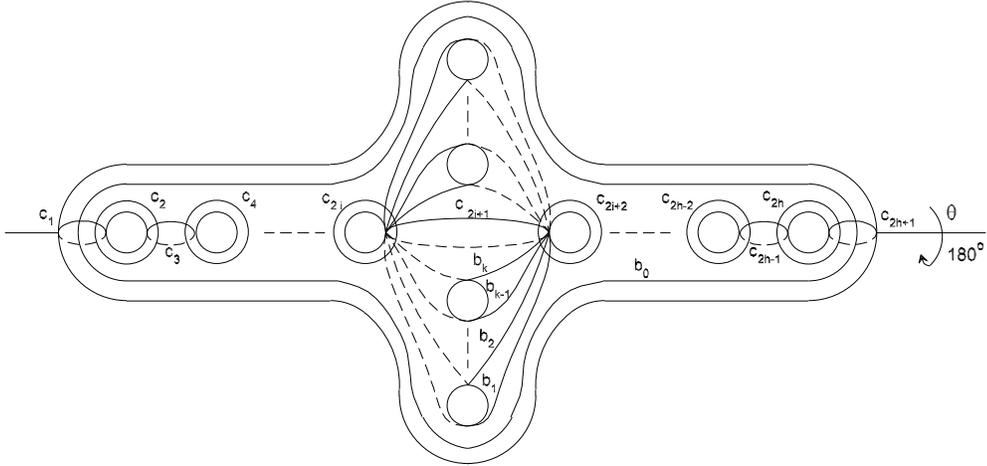,width=5.50in,clip=}
     \caption{The cycles used in the expression of $\theta$}
     \label{simplecase.fig}
 \end{figure}

Figure \ref{simplecase.fig} shows the cycles that are used in
expressing $\theta$ as a product of positive Dehn twists which is
stated in the next theorem.

\begin{thm} \label{simplecase.thm}
 The positive Dehn twist expression for the
involution $\theta$ on the surface $\Sigma_{h+k}$ shown in Figure
\ref{simplecase.fig} is given by
\[
\theta =c_{2i+2}\cdots c_{2h}c_{2h+1}c_{2i}\cdots
c_{2}c_{1}b_{0}c_{2h+1}c_{2h}\cdots c_{2i+2}c_{1}c_{2}\cdots
c_{2i}b_{1}b_{2}\cdots b_{k-1}b_{k}c_{2i+1}.
\]
\end{thm}

See \cite{Gu1} for the proof.\\

The order of the twists is from right to left, i.e., $c_{2i+1}$ is
applied first.

\section{Main Results}

\begin{figure}[htbp]
     \centering  \leavevmode
     \psfig{file=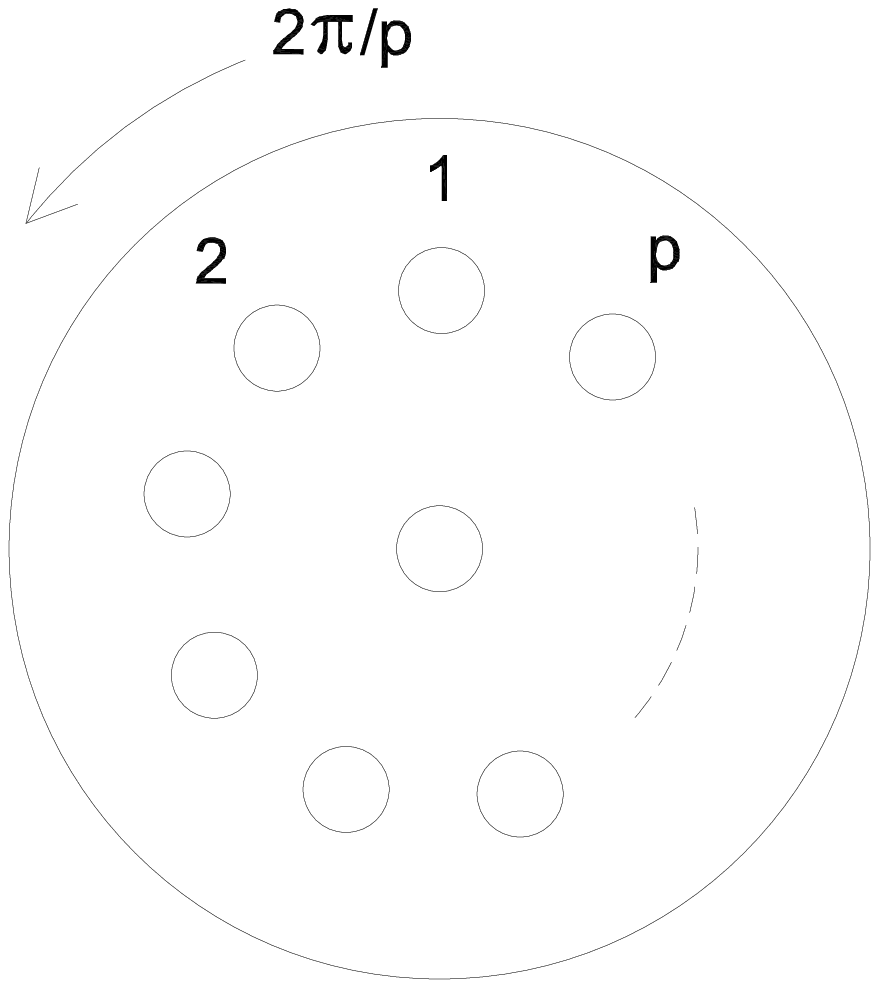,width=2.70in,clip=}
      \caption{$2\pi/p$ rotation $\phi_p$ on $\Sigma_{1+p}$ }
    \label{protation.fig}
 \end{figure}

It is very easy to define many of the elements of finite order in
the mapping class group qualitatively using geometry. In this
section we will consider those which are realized as rotations on
the surface $\Sigma_{1+p}$ through an angle of $2\pi /p$, Figure
\ref{protation.fig}.

It is a very difficult task in general,  however, to find a
positive Dehn twist expression for a given finite order element in
the mapping class group. In this section we will find explicit
positive Dehn twist expressions for rotations of the kind
described above. To achieve that, we will make use of the
expression for the involution $\theta$ that is given in Theorem
\ref{simplecase.thm}. $p$ represents a positive odd integer
throughout this section.

\begin{figure}[htbp]
     \centering  \leavevmode
     \psfig{file=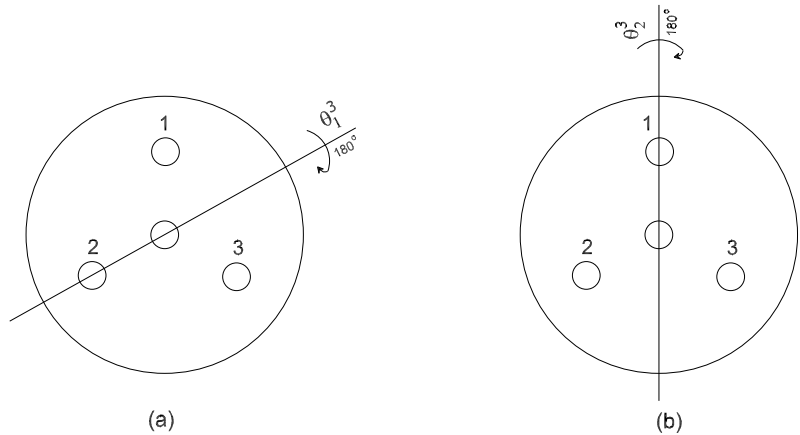,width=4.250in,clip=}
      \caption{Definition of the involutions $\theta_1^3$ and $\theta_2^3$ on $\Sigma_{1+3}$ }
    \label{theta1andtheta2together.fig}
 \end{figure}

Let $\phi_p$ denote the $2\pi /p$ rotation on the surface
$\Sigma_{1+p}$ as shown in Figure \ref{protation.fig}.

Let $\theta_1^p$ be the involution defined as $180^\circ $
rotation about the axis through the center of the surface
$\Sigma_{1+p}$ and the hole number $(p+1)/2,$ Figure
\ref{theta1andtheta2together.fig} (a). Let $\theta_2^p$ be the
involution defined as $180^\circ $ rotation about the axis through
the center of the surface $\Sigma_{1+p}$ and the hole number $1$,
Figure \ref{theta1andtheta2together.fig} (b).

\begin{figure}[htbp]
     \centering  \leavevmode
     \psfig{file=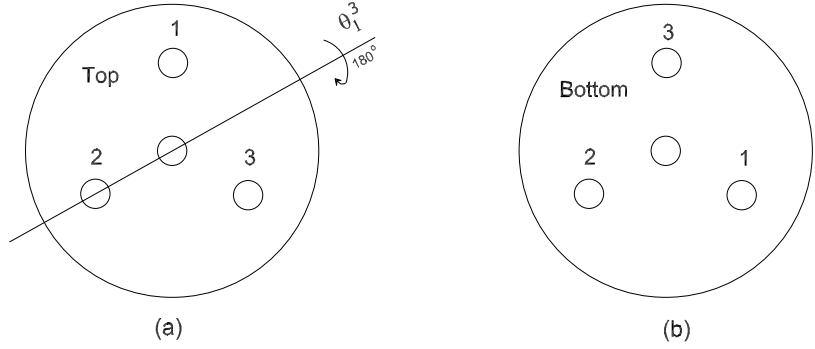,width=4.50in,clip=}
      \caption{Involution $\theta_1^3$ applied to $\Sigma_{1+3}$}
    \label{flip1.fig}
 \end{figure}

\begin{figure}[htbp]
     \centering  \leavevmode
     \psfig{file=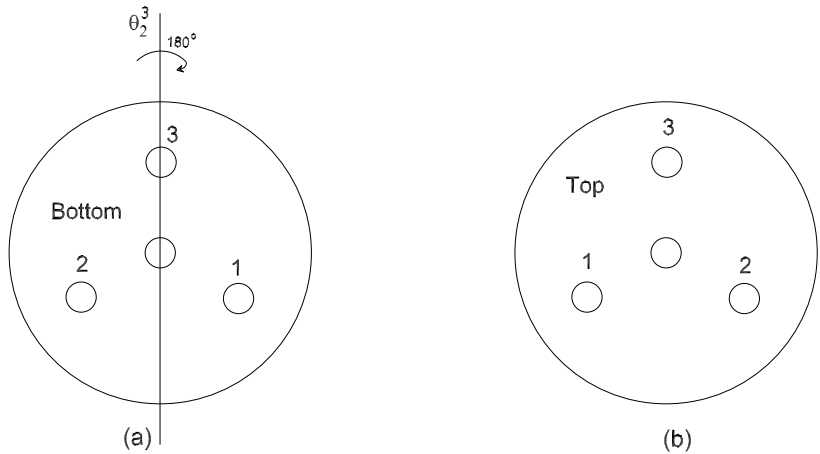,width=4.50in,clip=}
      \caption{Involution $\theta_2^3$ applied to the result from Figure \ref{flip1.fig}}
    \label{flip2.fig}
 \end{figure}

 It's not difficult to see that application of $\theta_1^p$ followed
 by $\theta_2^p$ results in $2\pi /p$  rotation about the center of the
 surface. Application of $\theta_1^p$ on the surface gives the
 result shown in Figure \ref{flip1.fig} (b), and $\theta_2^p$ takes
 this figure to the result shown in Figure \ref{flip2.fig} (b).

\begin{figure}[htbp]
     \centering  \leavevmode
     \psfig{file=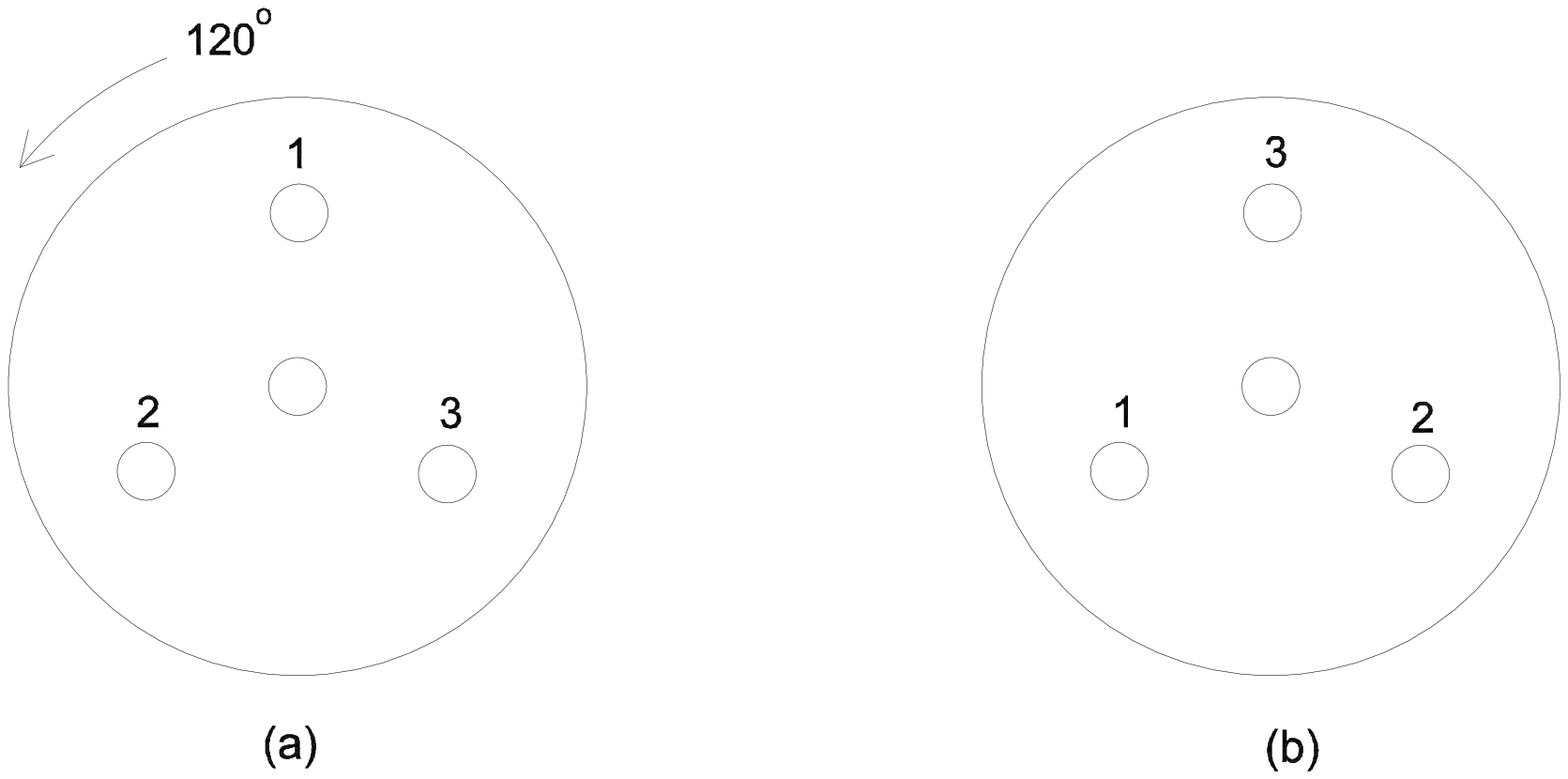,width=4.50in,clip=}
      \caption{$2\pi/3$ rotation $\phi_3$ on $\Sigma_{1+3}$ }
    \label{120rotation.fig}
 \end{figure}

Note that Figure \ref{flip2.fig} (b) and  Figure
\ref{120rotation.fig} (b) are the same, namely the result of
successive applications of the involutions $\theta_1^p$ and
$\theta_2^p$ is the same as the action of $\phi_p$ on the surface
$\Sigma_{1+p}$. From this we conclude that
\[\phi_p=\theta_2^p\theta_1^p\]
on the surface $\Sigma_{1+p}$.

\vspace{.5in}

Now the question is how to find a positive Dehn twist product for
$\phi_p$ . The answer  is simply juxtaposing the positive Dehn
twist expressions for $\theta_1^p$ and $\theta_2^p$. Therefore, we
first need to get those expressions for $\theta_1^p$ and
$\theta_2^p$ using Theorem \ref{simplecase.thm}.

\begin{figure}[htbp]
     \centering  \leavevmode
     \psfig{file=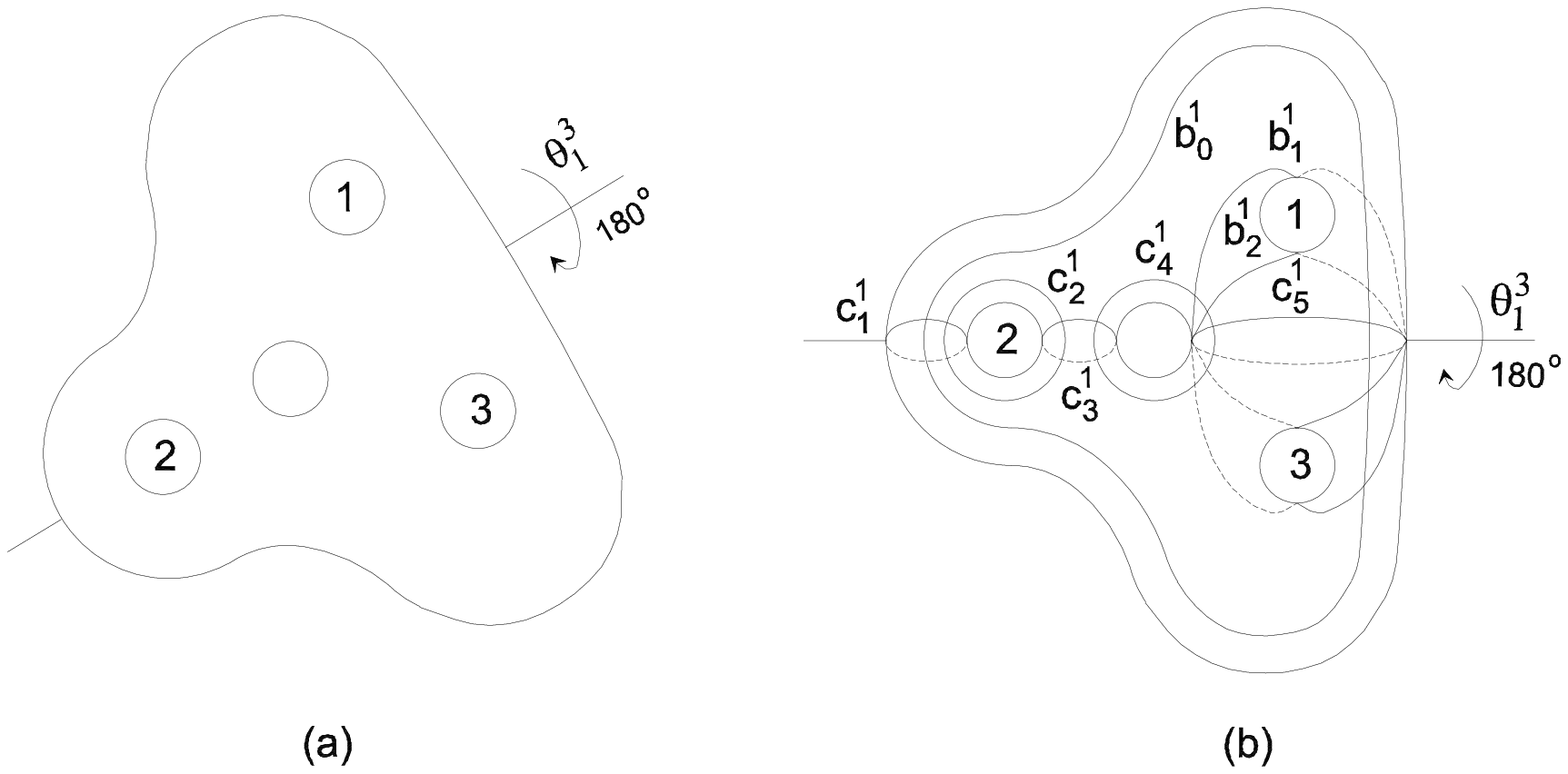,width=6.250in,clip=}
      \caption{Involution $\theta_1^3$ on $\Sigma_{1+3}$ and the cycles that realize it }
    \label{theta1cyclesnotflat.fig}
 \end{figure}

\begin{figure}[htbp]
     \centering  \leavevmode
     \psfig{file=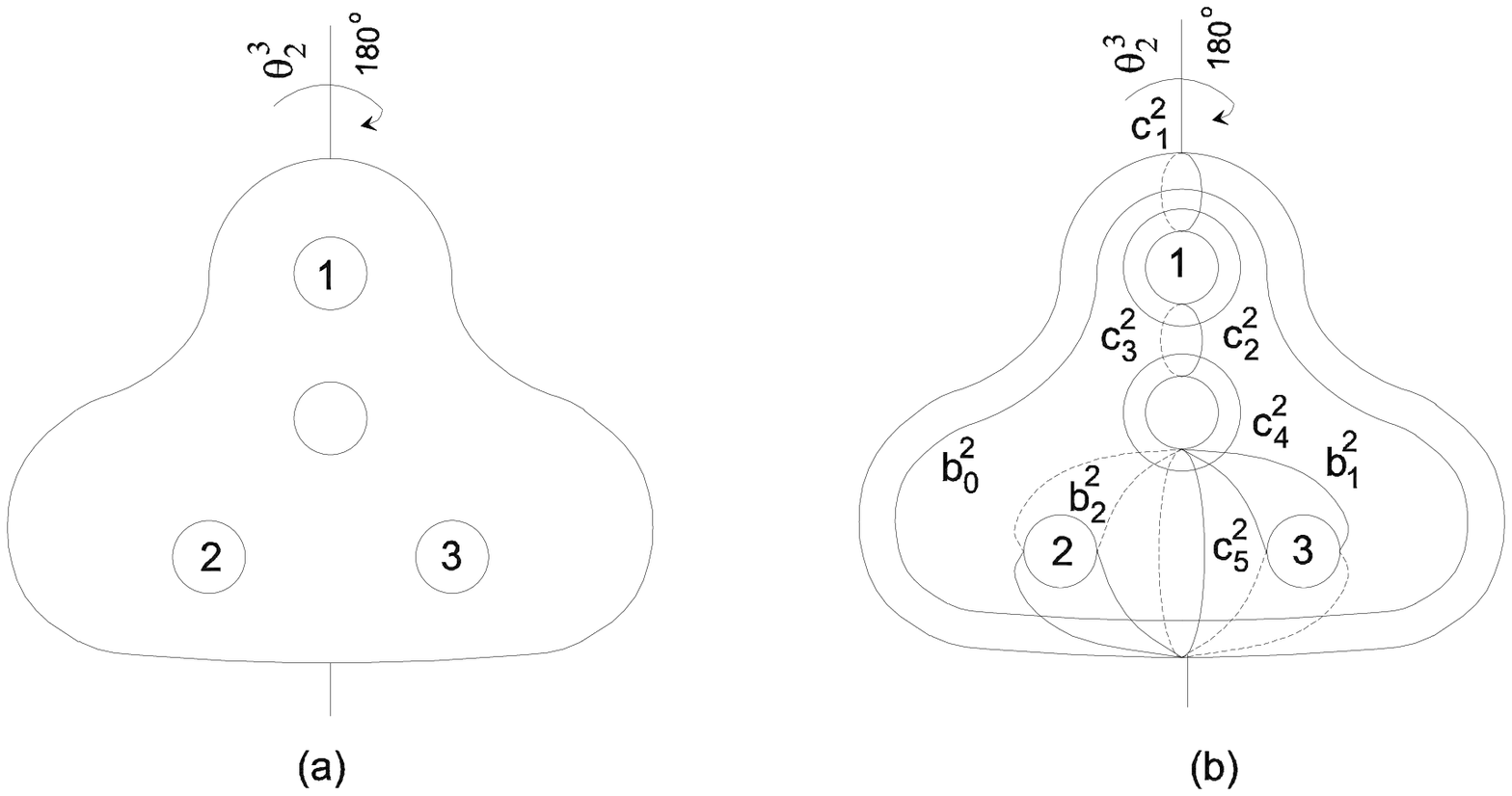,width=6.250in,clip=}
      \caption{Involution $\theta_2^3$ on $\Sigma_{1+3}$ and the cycles that realize it }
    \label{theta2cyclesnotflat.fig}
 \end{figure}

Figure \ref{theta1cyclesnotflat.fig} shows the cycles that realize
$\theta_1^3$ on the surface $\Sigma_{1+3}$ using Theorem
\ref{simplecase.thm}. Similarly, Figure
\ref{theta2cyclesnotflat.fig}  shows the cycles that realize
$\theta_2^3$ on the surface $\Sigma_{1+3}$ using the same theorem.

 According to Theorem
\ref{simplecase.thm} we have
\[\theta_1^3=c_{4}^{1}c_{3}^{1}c_{2}^{1}c_{1}^{1}b_{0}^{1}c_{1}^{1}c_{2}^{1}c_{3}^{1}c_{4}^{1}b_{1}^{1}b_{2}^{1}c_{5}^{1}
\]
and
\[\theta_2^3=c_{4}^{2}c_{3}^{2}c_{2}^{2}c_{1}^{2}b_{0}^{2}c_{1}^{2}c_{2}^{2}c_{3}^{2}c_{4}^{2}b_{1}^{2}b_{2}^{2}c_{5}^{2}.
\]
Therefore
\[\phi_3=\theta_2^3\theta_1^3=c_{4}^{2}c_{3}^{2}c_{2}^{2}c_{1}^{2}b_{0}^{2}c_{1}^{2}c_{2}^{2}c_{3}^{2}c_{4}^{2}b_{1}^{2}b_{2}^{2}c_{5}^{2}
c_{4}^{1}c_{3}^{1}c_{2}^{1}c_{1}^{1}b_{0}^{1}c_{1}^{1}c_{2}^{1}c_{3}^{1}c_{4}^{1}b_{1}^{1}b_{2}^{1}c_{5}^{1}.
\]
Since $\phi_3^3=1$, we have
\[
\left(
c_{4}^{2}c_{3}^{2}c_{2}^{2}c_{1}^{2}b_{0}^{2}c_{1}^{2}c_{2}^{2}c_{3}^{2}c_{4}^{2}b_{1}^{2}b_{2}^{2}c_{5}^{2}c_{4}^{1}c_{3}^{1}c_{2}^{1}c_{1}^{1}b_{0}^{1}c_{1}^{1}c_{2}^{1}c_{3}^{1}c_{4}^{1}b_{1}^{1}b_{2}^{1}c_{5}^{1}\right)
^{3}=1
\]
in $M_4$, the mapping class group of the surface $\Sigma_{4}$ .\\

\begin{figure}[htbp]
     \centering  \leavevmode
     \psfig{file=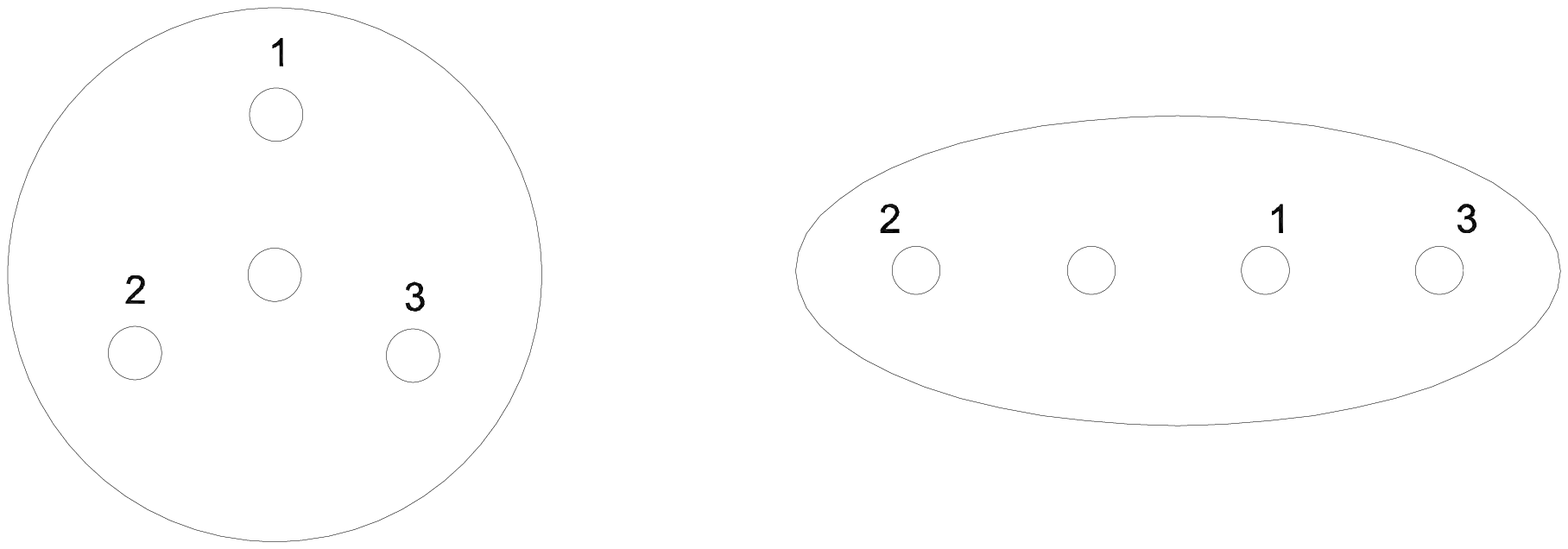,width=3.80in,clip=}
      \caption{Identification between the two surfaces $\Sigma_4$ and $\Sigma_{1+3}$ }
    \label{round2flat.fig}
 \end{figure}

 Next, we will demonstrate the cycles used in the expression for $\phi_3$ on a regular genus 4 surface
 $\Sigma_4$ and generalize it. For this, we need an identification between the two
 surfaces $\Sigma_4$ and $\Sigma_{1+3}$. Figure
 \ref{round2flat.fig} shows that identification.

\begin{figure}[htbp]
     \centering  \leavevmode
     \psfig{file=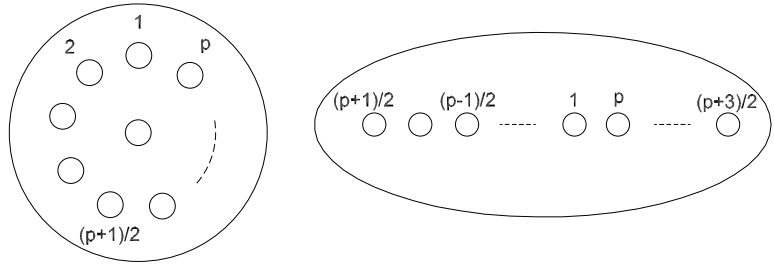,width=4.50in,clip=}
      \caption{Identification between the two surfaces for general case }
    \label{pround2flat.fig}
 \end{figure}

Figure \ref{pround2flat.fig} shows the same identification for
general case. The hole with no number on the right is the hole in
the center on the left.

\begin{figure}[htbp]
     \centering  \leavevmode
     \psfig{file=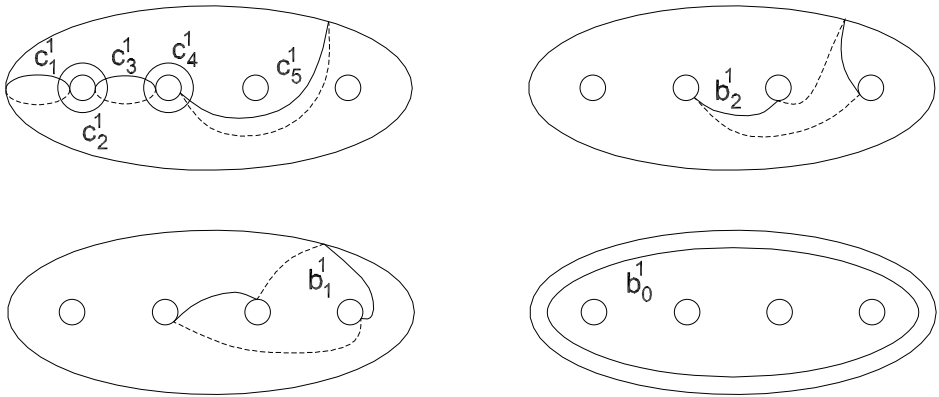,width=5.0in,clip=}
      \caption{The cycles realizing $\theta_1^3$ on $\Sigma_{4}$ }
    \label{theta1cyclesflat.fig}
 \end{figure}

\begin{figure}[htbp]
     \centering  \leavevmode
     \psfig{file=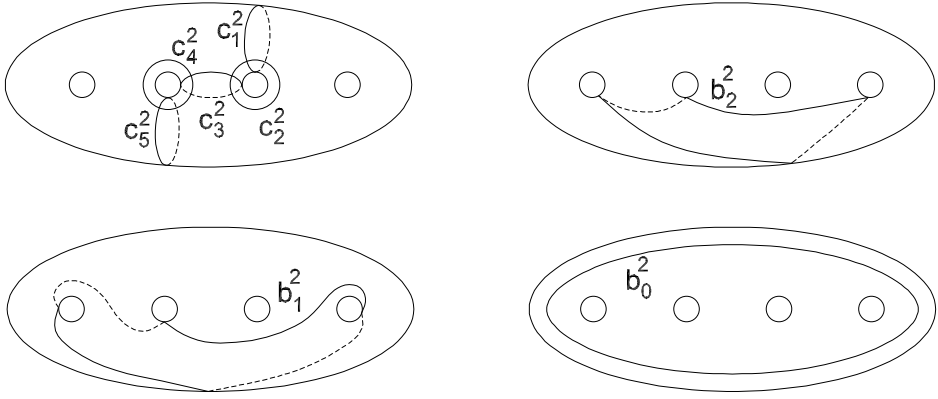,width=5.0in,clip=}
      \caption{The cycles realizing $\theta_2^3$ on $\Sigma_{4}$ }
    \label{theta2cyclesflat.fig}
 \end{figure}

Figure \ref{theta1cyclesflat.fig} shows the cycles in Figure
\ref{theta1cyclesnotflat.fig} (b) via the identification
established in Figure \ref{round2flat.fig}. Similarly, Figure
\ref{theta2cyclesflat.fig} shows the cycles in Figure
\ref{theta2cyclesnotflat.fig} (b) via the identification
established in Figure \ref{round2flat.fig}.

In general, the cycles realizing $\phi_p$ on the round surface
$\Sigma_{1+p}$ are carried onto the regular surface $\Sigma_{p+1}$
using the
identification in Figure \ref{pround2flat.fig}.\\

\section{Applications}

Let $\Sigma_{g}$ be the $2$-dimensional, closed, compact, oriented
surface of genus $g>0$ and $M_g$ be its mapping class group, the
group of isotopy classes of all orientation-preserving
diffeomorphisms $\Sigma_{g}\rightarrow \Sigma_{g}.$ It is a
well-known fact that any word in the mapping class group $M_g$
that contains positive exponents only and is equal to the identity
element defines a symplectic Lefschetz fibration $
X^{4}\rightarrow S^{2}$.  In this section we compute the
homeomorphism invariants of the  symplectic Lefschetz fibrations
that are defined by the words $\phi_p^p=1$, $p-$ odd, in
$M_{p+1}$.

\begin{figure}[htbp]
     \centering  \leavevmode
     \psfig{file=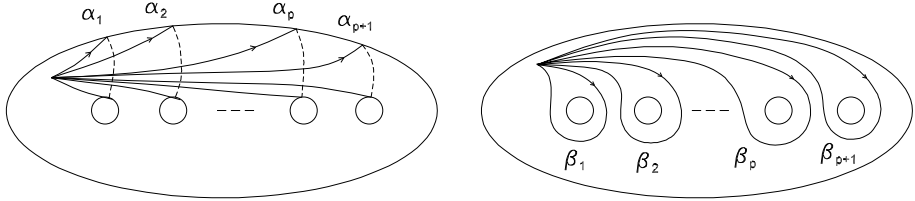,width=5.0in,clip=}
      \caption{Standard generators of $\pi_1\left(\Sigma_{p+1}\right)$ }
    \label{fundgroupgenerators.fig}
 \end{figure}

 First, we will show that the $4-$ manifold $X$ carrying the
 symplectic Lefschetz fibration structure $X^{4}\rightarrow S^{2}$
is simply connected. A well known fact in theory of Lefschetz
fibrations states that $\pi_1(X)$ is isomorphic  to the quotient
of  $\pi_1\left(\Sigma_{p+1}\right)$ by the normal subgroup
generated by the vanishing cycles, the cycles about which positive
Dehn twists in the expression $\phi_p^p=1$ are performed. The
vanishing cycles are seen as elements of
$\pi_1\left(\Sigma_{p+1}\right)$ in this section, not as elements
of $M_{g+1}$.

\begin{figure}[htbp]
     \centering  \leavevmode
     \psfig{file=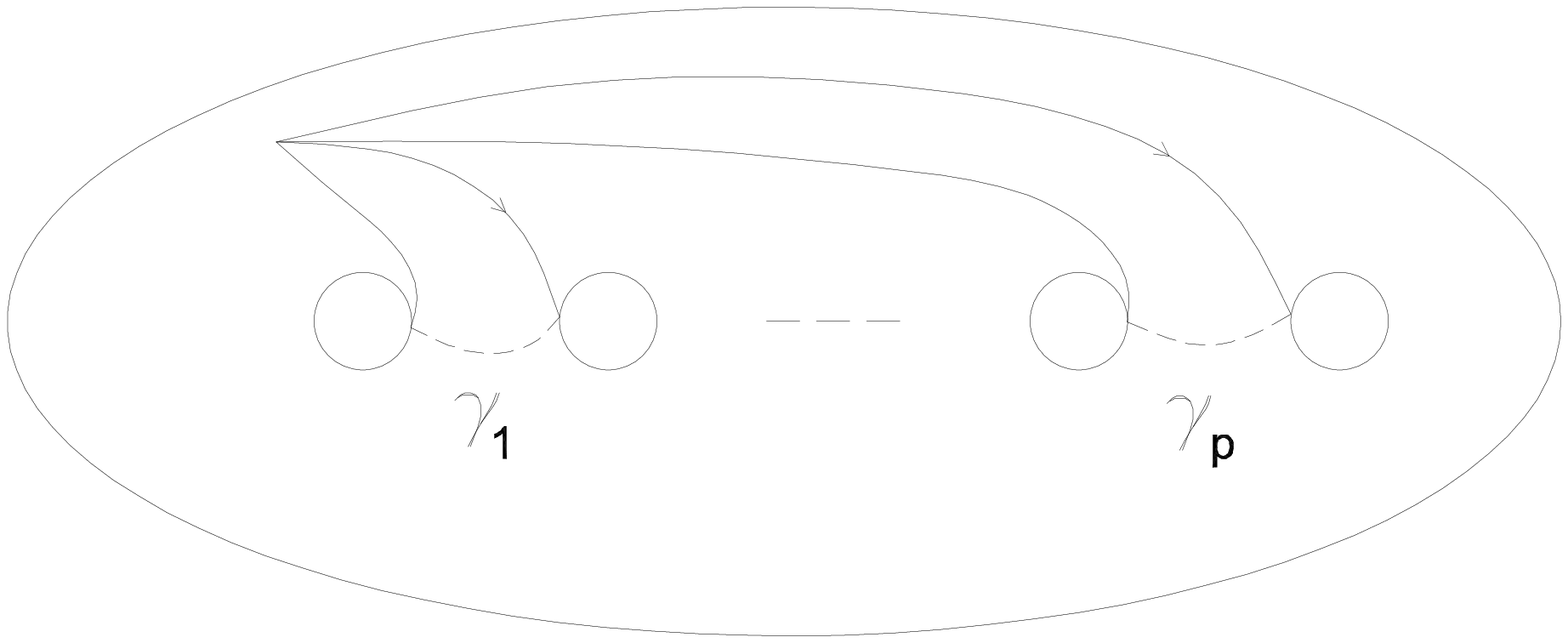,width=2.50in,clip=}
      \caption{$\gamma_{i}=\alpha_{i}\alpha_{i+1}^{-1},\ 1\leq i \leq p$ }
    \label{gammacycles.fig}
 \end{figure}

Let $\{\alpha_{1},\beta_{1},\cdots, \beta_{p+1},\alpha_{p+1}\}$ be
the set of standard $2p+2$  generators of
$\pi_1\left(\Sigma_{p+1}\right)$ as shown in  Figure
\ref{fundgroupgenerators.fig}.  We will show that the subgroup of
$\pi_1\left(\Sigma_{p+1}\right)$ generated by the vanishing cycles
includes all the generators of $\pi_1\left(\Sigma_{p+1}\right)$ by
showing that each generator is equal to 1 in the quotient group.
We will  use the elements in $\pi_1\left(\Sigma_{p+1}\right)$ that
are  shown in Figure \ref{gammacycles.fig} along with the standard
generators to express the vanishing cycles as elements of
$\pi_1\left(\Sigma_{p+1}\right)$. The expressions should be read
from right to left.

If $p=3$ then

\[\theta_1^3=c_{4}^{1}c_{3}^{1}c_{2}^{1}c_{1}^{1}b_{0}^{1}c_{1}^{1}c_{2}^{1}c_{3}^{1}c_{4}^{1}b_{1}^{1}b_{2}^{1}c_{5}^{1}
\]
\[\theta_2^3=c_{4}^{2}c_{3}^{2}c_{2}^{2}c_{1}^{2}b_{0}^{2}c_{1}^{2}c_{2}^{2}c_{3}^{2}c_{4}^{2}b_{1}^{2}b_{2}^{2}c_{5}^{2}
\]
and  we have the following expressions in
$\pi_1\left(\Sigma_{4}\right)$ for the vanishing cycles shown in
Figures \ref{theta1cyclesflat.fig} and \ref{theta2cyclesflat.fig}:

 \[c_{1}^{1}=\alpha_{1},c_{2}^{1}=\beta_{1},c_{3}^{1}=\gamma_{1},c_{4}^{1}=\beta_{2},
 c_{5}^{1}=\beta_{3}\alpha_{3}^{-1}\beta_{3}^{-1}\gamma_{2}^{-1}\]
\[b_{0}^{1}=\beta_{4}^{-1}\beta_{3}^{-1} \beta_{2}^{-1}\beta_{1}^{-1}\]
\[b_{1}^{1}=\alpha_{3}\beta_{4}^{-1}\gamma_{3}^{-1}\beta_{3}^{-1}\gamma_{2}^{-1}\]
\[b_{2}^{1}=\beta_{3}\alpha_{3}\gamma_{3}^{-1}\beta_{3}^{-1}\gamma_{2}^{-1}\]
\[c_{1}^{2}=\alpha_{3},c_{2}^{2}=\beta_{3},c_{3}^{2}=\gamma_{2},c_{4}^{2}=\beta_{2},c_{5}^{2}
=\gamma_{1}^{-1}\beta_{1}^{-1}\alpha_{1}\beta_{1}\]
\[b_{0}^{2}=\beta_{4}^{-1}\beta_{3}^{-1} \beta_{2}^{-1}\beta_{1}^{-1}\]
\[b_{1}^{2}=\gamma_{1}\beta_{2}\beta_{3}
\beta_{4}^{-1}\gamma_{3}^{-1}\beta_{3}^{-1}\gamma_{2}^{-1}\beta_{2}^{-1}
\gamma_{1}^{-1}\beta_{1}^{-1}\alpha_{1}\]
\[b_{2}^{2}=\gamma_{1}\beta_{2}\beta_{3}
\gamma_{3}^{-1}\beta_{3}^{-1}\gamma_{2}^{-1}\beta_{2}^{-1}
\gamma_{1}^{-1}\beta_{1}^{-1}\alpha_{1}\beta_{1}\]\\

 Let $q=\frac{p+1}{2},$
i.e., $p=2q-1$. Then for $p\geq 5$, and hence $q\geq 3,$

\begin{eqnarray}    \label{theta p}
\theta_2^p=
c_{4}^{2}c_{3}^{2}c_{2}^{2}c_{1}^{2}b_{0}^{2}c_{1}^{2}c_{2}^{2}c_{3}^{2}c_{4}^{2}b_{1}^{2}
\cdots b_{p-1}^{2}c_{5}^{2}
\end{eqnarray}
\begin{eqnarray*}
\theta_1^p=
c_{4}^{1}c_{3}^{1}c_{2}^{1}c_{1}^{1}b_{0}^{1}c_{1}^{1}c_{2}^{1}c_{3}^{1}c_{4}^{1}b_{1}^{1}
\cdots b_{p-1}^{1}c_{5}^{1}
\end{eqnarray*}
 and we have the following expressions in
$\pi_1\left(\Sigma_{p+1}\right)$ for the vanishing cycles: \\

 $\theta_1^p$ cycles:

 \[c_{1}^{1}=\alpha_{1},c_{2}^{1}=\beta_{1},c_{3}^{1}=\gamma_{1},c_{4}^{1}=\beta_{2}\]

\[c_{5}^{1}=\beta_{3}\cdots
\beta_{q}\beta_{q+1}\alpha_{q+1}^{-1}\beta_{q+1}^{-1}\gamma_{q}^{-1}\cdots
\beta_{4}^{-1}\gamma_{3}^{-1}\beta_{3}^{-1}\gamma_{2}^{-1}\]

\[b_{0}^{1}=\beta_{p+1}^{-1}\beta_{p}^{-1}\cdots \beta_{2}^{-1}\beta_{1}^{-1}\]

\[b_{1}^{1}=\alpha_{3}\beta_{p+1}^{-1}\gamma_{p}^{-1}\beta_{p}^{-1}\cdots
\beta_{4}^{-1}\gamma_{3}^{-1}\beta_{3}^{-1}\gamma_{2}^{-1}\]

\[b_{2m}^{1}=\beta_{3}\cdots \beta_{m+1}
\beta_{m+2}\alpha_{m+2}\gamma_{p-m+1}^{-1}\beta_{p-m+1}^{-1}\gamma_{p-m}^{-1}\cdots
\beta_{4}^{-1}\gamma_{3}^{-1}\beta_{3}^{-1}\gamma_{2}^{-1},q-1\geq
m\geq 1\]

\[b_{2m-1}^{1}=\beta_{3}\cdots \beta_{m}
\beta_{m+1}\alpha_{m+2}\beta_{p-m+2}^{-1}\gamma_{p-m+1}^{-1}\beta_{p-m+1}^{-1}\cdots
\beta_{4}^{-1}\gamma_{3}^{-1}\beta_{3}^{-1}\gamma_{2}^{-1},q-1\geq m\geq 2\] \\

 $\theta_2^p$ cycles:

\[c_{1}^{2}=\alpha_{q+1},c_{2}^{2}=\beta_{q+1},c_{4}^{2}=\beta_{2},c_{5}^{2}=\gamma_{1}^{-1}\beta_{1}^{-1}\alpha_{1}\beta_{1}\]

\[c_{3}^{2}=\beta_{3}\cdots\beta_{q-1} \beta_{q}\gamma_{q}^{-1}\beta_{q}^{-1}\cdots
 \gamma_{3}^{-1}\beta_{3}^{-1}\gamma_{2}^{-1}\]

\[b_{0}^{2}=\beta_{p+1}^{-1}\beta_{p}^{-1}\cdots \beta_{2}^{-1}\beta_{1}^{-1}\]

\[b_{p-1}^{2}=\gamma_{1}\beta_{2}\beta_{3}\cdots
\beta_{p-1}\beta_{p}\gamma_{p}^{-1}\beta_{p}^{-1}\gamma_{p-1}^{-1}\cdots
\gamma_{1}^{-1}\beta_{1}^{-1}\alpha_{1}\beta_{1}\]

\[b_{p-2}^{2}=\gamma_{1}\beta_{2}\beta_{3}\cdots
\beta_{p-1}\beta_{p}\beta_{p+1}^{-1}\gamma_{p}^{-1}\beta_{p}^{-1}\gamma_{p-1}^{-1}\cdots
\gamma_{1}^{-1}\beta_{1}^{-1}\alpha_{1}\]

\[b_{2m-1}^{2}=\beta_{1+q-m}^{-1}\gamma_{q-m}^{-1}\cdots
\beta_{4}^{-1}\gamma_{3}^{-1}\beta_{3}^{-1}\gamma_{2}^{-1}\beta_{2}\beta_{3}\cdots
\beta_{q+m}\gamma_{q+m+1}\cdots
\gamma_{p-1}\gamma_{p}\beta_{p+1}^{-1}\gamma_{p}^{-1}\beta_{p}^{-1}\gamma_{p-1}^{-1}\cdots
\gamma_{1}^{-1}\beta_{1}^{-1}\alpha_{1},\] $q-2\geq m\geq 1$

\[b_{2m}^{2}=\gamma_{q-m}^{-1}\beta_{q-m}^{-1}\cdots
\gamma_{3}^{-1}\beta_{3}^{-1}\gamma_{2}^{-1}\beta_{2}\beta_{3}\cdots
\beta_{q+m+1}\gamma_{q+m+1}\cdots
\gamma_{p-1}\gamma_{p}\beta_{p+1}^{-1}\gamma_{p}^{-1}\beta_{p}^{-1}\gamma_{p-1}^{-1}\cdots
\gamma_{1}^{-1}\beta_{1}^{-1}\alpha_{1},\] $ q-2\geq m\geq 1.$\\

All of the vanishing cycles above are equal to 1 in the quotient
$\pi_1\left(\Sigma_{p+1}\right)/N,$ where $N$ is the normal
subgroup generated by the vanishing cycles. Therefore we can
cancel some of the standard generators right away:
$\alpha_{1},\beta_{1},\beta_{2},\alpha_{q+1},\beta_{q+1}$ can be
cancelled  because they are equal to the cycles
$c_{1}^{1},c_{2}^{1},c_{4}^{1},c_{1}^{2},c_{2}^{2}$, respectively.
Also, $\alpha_{2}$ can be cancelled because
$c_{3}^{1}=\gamma_{1}=\alpha_{1}\alpha_{2}^{-1}$. There remains
$2p-4$ elements, which are

\begin{equation} \label{remaining generators.list}
\alpha_{3},\alpha_{4},\cdots
,\alpha_{q},\alpha_{q+2},\cdots ,\alpha_{p},\alpha_{p+1} \\
\beta_{3},\beta_{4},\cdots ,\beta_{q},\beta_{q+2}\cdots
,\beta_{p},\beta_{p+1}.
\end{equation}\\Now,

\[b_{1}^{1}=\alpha_{3}\beta_{p+1}^{-1}\gamma_{p}^{-1}\beta_{p}^{-1}\cdots
\beta_{4}^{-1}\gamma_{3}^{-1}\beta_{3}^{-1}\gamma_{2}^{-1}=1\]\\
and

\[b_{2}^{1}=\beta_{3}\alpha_{3}\gamma_{p}^{-1}\beta_{p}^{-1}\cdots
\beta_{4}^{-1}\gamma_{3}^{-1}\beta_{3}^{-1}\gamma_{2}^{-1}=1\]\\
give

\[\alpha_{3}\beta_{p+1}^{-1}=\beta_{3}\alpha_{3},\ \text{i.e.},\
\beta_{3}=\alpha_{3}\beta_{p+1}^{-1}\alpha_{3}^{-1}.\]\\
Similarly,

\[b_{3}^{1}=\beta_{3}\alpha_{4}\beta_{p}^{-1}\gamma_{p-1}^{-1}\beta_{p-1}^{-1}\cdots
\beta_{4}^{-1}\gamma_{3}^{-1}\beta_{3}^{-1}\gamma_{2}^{-1}=1\]
\\ and

\[b_{4}^{1}=\beta_{3}\beta_{4}\alpha_{4}\gamma_{p-1}^{-1}\beta_{p-1}^{-1}\gamma_{p-2}^{-1}\cdots
\beta_{4}^{-1}\gamma_{3}^{-1}\beta_{3}^{-1}\gamma_{2}^{-1}=1\]\\
give

\[\beta_{3}\beta_{4}\alpha_{4}=\beta_{3}\alpha_{4}\beta_{p}^{-1},\ \text{i.e.},\ \beta_{4}=\alpha_{4}\beta_{p}^{-1}\alpha_{4}^{-1}.
\]\\
In general

\[b_{2m-1}^{1}=\beta_{3}\cdots
\beta_{m+1}\alpha_{m+2}\beta_{p-m+2}^{-1}\gamma_{p-m+1}^{-1}\cdots
\beta_{4}^{-1}\gamma_{3}^{-1}\beta_{3}^{-1}\gamma_{2}^{-1}=1\]
\\ and

\[b_{2m}^{1}=\beta_{3}\cdots
\beta_{m+2}\alpha_{m+2}\gamma_{p-m+1}^{-1}\beta_{p-m+1}^{-1}\gamma_{p-m}^{-1}\cdots
\beta_{4}^{-1}\gamma_{3}^{-1}\beta_{3}^{-1}\gamma_{2}^{-1}= 1\]\\
give

\begin{equation} \label{b1-conjugation.relation}
 \beta_{m+2}=\alpha_{m+2}\beta_{p-m+2}^{-1}\alpha_{m+2}^{-1},\ q-1\geq m\geq
 1.
 \end{equation}\\

A similar conjugation relation can be obtained by looking at the
cycles of $\theta_2^p$ in pairs as well. The first pair is
$b_{1}^{2}$ and $b_{2}^{2}$: \\

\[b_{1}^{2}=\beta_{q}^{-1}\gamma_{q-1}^{-1}\cdots
\gamma_{3}^{-1}\beta_{3}^{-1}\gamma_{2}^{-1}\beta_{2}\cdots
\beta_{q+1}\gamma_{q+2}\cdots
\gamma_{p}\beta_{p+1}^{-1}\gamma_{p}^{-1}\beta_{p}^{-1}\gamma_{p-1}^{-1}\cdots
\gamma_{1}^{-1}\beta_{1}^{-1}\alpha_{1}=1\]\\ and

\[b_{2}^{2}=\gamma_{q-1}^{-1}\beta_{q-1}^{-1}\cdots
\gamma_{3}^{-1}\beta_{3}^{-1}\gamma_{2}^{-1}\beta_{2}\cdots
\beta_{q+2}\gamma_{q+2}\cdots
\gamma_{p}\beta_{p+1}^{-1}\gamma_{p}^{-1}\beta_{p}^{-1}\gamma_{p-1}^{-1}\cdots
\gamma_{1}^{-1}\beta_{1}^{-1}\alpha_{1}=1\]\\ result in

\[\beta_{q}^{-1}\gamma_{q-1}^{-1}\cdots
\gamma_{3}^{-1}\beta_{3}^{-1}\gamma_{2}^{-1}\beta_{2}\cdots
\beta_{q+1}=\gamma_{q-1}^{-1}\beta_{q-1}^{-1}\cdots
\gamma_{3}^{-1}\beta_{3}^{-1}\gamma_{2}^{-1}\beta_{2}\cdots
\beta_{q+1}\beta_{q+2},\]\\ which can be written as

\[\left(\gamma_{q-1}^{-1}\beta_{q-1}^{-1}\cdots
\gamma_{3}^{-1}\beta_{3}^{-1}\gamma_{2}^{-1}\beta_{2}\cdots
\beta_{q+1}\right)^{-1}\beta_{q}^{-1}\left(\gamma_{q-1}^{-1}\cdots
\gamma_{3}^{-1}\beta_{3}^{-1}\gamma_{2}^{-1}\beta_{2}\cdots
\beta_{q+1}\right)=\beta_{q+2}.\]\\
Similarly, from

\[b_{3}^{2}=\beta_{q-1}^{-1}\gamma_{q-2}^{-1}\cdots
\gamma_{3}^{-1}\beta_{3}^{-1}\gamma_{2}^{-1}\beta_{2}\cdots
\beta_{q+2}\gamma_{q+3}\cdots
\gamma_{p}\beta_{p+1}^{-1}\gamma_{p}^{-1}\beta_{p}^{-1}\gamma_{p-1}^{-1}\cdots
\gamma_{1}^{-1}\beta_{1}^{-1}\alpha_{1}=1\]\\ and

\[b_{4}^{2}=\gamma_{q-2}^{-1}\beta_{q-2}^{-1}\cdots
\gamma_{3}^{-1}\beta_{3}^{-1}\gamma_{2}^{-1}\beta_{2}\cdots
\beta_{q+3}\gamma_{q+3}\cdots
\gamma_{p}\beta_{p+1}^{-1}\gamma_{p}^{-1}\beta_{p}^{-1}\gamma_{p-1}^{-1}\cdots
\gamma_{1}^{-1}\beta_{1}^{-1}\alpha_{1}=1\]\\ we obtain

\[\beta_{q-1}^{-1}\gamma_{q-2}^{-1}\cdots
\gamma_{3}^{-1}\beta_{3}^{-1}\gamma_{2}^{-1}\beta_{2}\cdots
\beta_{q+2}=\gamma_{q-2}^{-1}\beta_{q-2}^{-1}\cdots
\gamma_{3}^{-1}\beta_{3}^{-1}\gamma_{2}^{-1}\beta_{2}\cdots
\beta_{q+2}\beta_{q+3},\]\\ which can be written as

\[\left(\gamma_{q-2}^{-1}\beta_{q-2}^{-1}\cdots
\gamma_{3}^{-1}\beta_{3}^{-1}\gamma_{2}^{-1}\beta_{2}\cdots\beta_{q+2}\right)^{-1}
\beta_{q-1}^{-1}\left(\gamma_{q-2}^{-1}\cdots
\gamma_{3}^{-1}\beta_{3}^{-1}\gamma_{2}^{-1}\beta_{2}\cdots
\beta_{q+2}\right)= \beta_{q+3}.\]\\
In general, from

\[b_{2m-1}^{2}=\beta_{1+q-m}^{-1}\gamma_{q-m}^{-1}\cdots
\gamma_{3}^{-1}\beta_{3}^{-1}\gamma_{2}^{-1}\beta_{2}\cdots
\beta_{q+m}\gamma_{q+m+1}\cdots
\gamma_{p}\beta_{p+1}^{-1}\gamma_{p}^{-1}\beta_{p}^{-1}\gamma_{p-1}^{-1}\cdots
\gamma_{1}^{-1}\beta_{1}^{-1}\alpha_{1}=1\]\\ and

\[b_{2m}^{2}=\gamma_{q-m}^{-1}\beta_{q-m}^{-1}\cdots
\gamma_{3}^{-1}\beta_{3}^{-1}\gamma_{2}^{-1}\beta_{2}\cdots
\beta_{q+m+1}\gamma_{q+m+1}\cdots
\gamma_{p}\beta_{p+1}^{-1}\gamma_{p}^{-1}\beta_{p}^{-1}\gamma_{p-1}^{-1}\cdots
\gamma_{1}^{-1}\beta_{1}^{-1}\alpha_{1}=1\]\\ we obtain

\[\beta_{1+q-m}^{-1}\gamma_{q-m}^{-1}\cdots
\gamma_{3}^{-1}\beta_{3}^{-1}\gamma_{2}^{-1}\beta_{2}\cdots
\beta_{q+m}=\gamma_{q-m}^{-1}\beta_{q-m}^{-1}\cdots
\gamma_{3}^{-1}\beta_{3}^{-1}\gamma_{2}^{-1}\beta_{2}\cdots
\beta_{q+m}\beta_{q+m+1}\]\\ which can be written as

\begin{equation} \label{b2-conjugation.relation}
\beta_{1+q-m}=\left(\gamma_{q-m}^{-1}\cdots
\gamma_{3}^{-1}\beta_{3}^{-1}\gamma_{2}^{-1}\beta_{2}\cdots
\beta_{q+m}\right)\beta_{q+m+1}^{-1}\left(\gamma_{q-m}^{-1}\beta_{q-m}^{-1}\cdots
\gamma_{3}^{-1}\beta_{3}^{-1}\gamma_{2}^{-1}\beta_{2}\cdots
\beta_{q+m}\right)^{-1},\\
\end{equation}
$ q-2\geq m\geq 1$\\

 Setting $m=1$ in \ref{b1-conjugation.relation}
we see that $\beta_{3}$ and $\beta_{p+1}^{-1}$ are conjugate. On
the other hand substituting $m=q-2$ in
\ref{b2-conjugation.relation} we see that $\beta_{3}$ and
$\beta_{p}^{-1}$ are conjugate. Likewise setting $m=2$ in
\ref{b1-conjugation.relation} and $m=q-3$ in
\ref{b2-conjugation.relation} we see that
$\beta_{4},\beta_{p}^{-1}$ and $\beta_{p-1}^{-1}$ are conjugate to
each other. Continuing this way we get the following list of
triples.

\begin{eqnarray*}
  \beta_{3} & \beta_{p+1}^{-1} & \beta_{p}^{-1} \\
  \beta_{4} & \beta_{p}^{-1} & \beta_{p-1}^{-1} \\
  \beta_{5} & \beta_{p-1}^{-1} & \beta_{p-2}^{-1} \\
  \cdots & \cdots & \cdots \\
  \beta_{q} & \beta_{q+3}^{-1} & \beta_{q+2}^{-1} \\
\end{eqnarray*}

Each row in this list contains elements that are conjugate to each
other. It is easy now to see that $\beta_{3},\beta_{4},
\cdots,\beta_{q},\beta_{q+2}^{-1},\beta_{q+3}^{-1},\cdots ,
\beta_{p}^{-1},\beta_{p+1}^{-1}$ are all conjugate to each other.
Finally, setting $m= q-1$ in \ref{b1-conjugation.relation} we see
that $\beta_{q+1}$ and $\beta_{q+2}^{-1}$ are conjugate. Since
$c_{2}^{2}=\beta_{q+1}=1$ we conclude that the second half of the
generators in \ref{remaining generators.list} are all identity in
the quotient. Instead of $\beta_{q+1}=1$ here we could have used
$\beta_{p+1}=1$, which is obtained from the pair of cycles
$b_{p-1}^{2},b_{p-2}^{2}$:

\[b_{p-1}^{2}=\gamma_{1}\beta_{2}\cdots
\beta_{p-1}\beta_{p}\gamma_{p}^{-1}\beta_{p}^{-1}\gamma_{p-1}^{-1}\cdots
\gamma_{1}^{-1}\beta_{1}^{-1}\alpha_{1}\beta_{1}=1\] and

\[b_{p-2}^{2}=\gamma_{1}\beta_{2}\cdots
\beta_{p-1}\beta_{p}\beta_{p+1}^{-1}\gamma_{p}^{-1}\beta_{p}^{-1}\gamma_{p-1}^{-1}\cdots
\gamma_{1}^{-1}\beta_{1}^{-1}\alpha_{1}=1\]\\ give

\[\gamma_{1}\beta_{2}\cdots \beta_{p-1}\beta_{p}=\gamma_{1}\beta_{2}\cdots
\beta_{p-1}\beta_{p}\beta_{p+1}^{-1},\ \text{i.e.,}\
\beta_{p+1}=1\]\\
using the fact that $c_{2}^{1}=\beta_{1}=1$.\\

 Therefore, the list of
generators reduces to

\[
\alpha_{3},\alpha_{4},\cdots ,\alpha_{q},\alpha_{q+2},\cdots
,\alpha_{p},\alpha_{p+1}.
\]\\
Before we proceed further let's recall that
$\gamma_{i}=\alpha_{i}\alpha_{i+1}^{-1}$. Therefore

\begin{eqnarray*}
\gamma _{i}\cdots \gamma _{j}=\alpha _{i}\alpha _{i+1}^{-1}\alpha
_{i+1}\alpha _{i+2}^{-1}\cdots \alpha _{j-1}\alpha _{j}^{-1}\alpha
_{j}\alpha _{j+1}^{-1}=\alpha _{i}\alpha
_{j+1}^{-1}\\
\gamma_{j}^{-1}\cdots
\gamma_{i}^{-1}=\alpha_{j+1}\alpha_{j}^{-1}\alpha_{j}\alpha_{j-1}^{-1}\cdots
\alpha_{i+2}\alpha_{i+1}^{-1}\alpha_{i+1}\alpha_{i}^{-1}=\alpha_{j+1}\alpha_{i}^{-1}
\end{eqnarray*}

 for $j>i$. Now,

 \[b_{1}^{1}=\alpha_{3}\beta_{p+1}^{-1}\gamma_{p}^{-1}\cdots
\beta_{4}^{-1}\gamma_{3}^{-1}\beta_{3}^{-1}\gamma_{2}^{-1}=\alpha_{3}\gamma_{p}^{-1}\cdots
\gamma_{3}^{-1}\gamma_{2}^{-1}=\alpha_{3}\alpha_{p+1}\alpha_{2}^{-1}=\alpha_{3}\alpha_{p+1}=1.\]\\
In general we have

\[b_{2m-1}^{1}=\beta _{3}\cdots \beta _{m+1}\alpha _{m+2}\beta
_{p-m+2}^{-1}\gamma _{p-m+1}^{-1}\cdots \beta _{4}^{-1}\gamma
_{3}^{-1}\beta _{3}^{-1}\gamma _{2}^{-1}=\alpha _{m+2}\gamma
_{p-m+1}^{-1}\cdots \gamma _{3}^{-1}\gamma _{2}^{-1}\]
\begin{equation} \label{alpha relations1}
 =\alpha _{m+2}\alpha _{p-m+2}\alpha
_{2}^{-1}=\alpha _{m+2}\alpha _{p-m+2}=1,\ q-1\geq m\geq 2.
\end{equation}

A similar list of products can be obtained using the cycles of
$\theta_{2}^{p}$ as follows:

\[
b_{1}^{2}=\beta_{q}^{-1}\gamma_{q-1}^{-1}\cdots
\gamma_{3}^{-1}\beta_{3}^{-1}\gamma_{2}^{-1}\beta_{2}\cdots
\beta_{q+1}\gamma_{q+2}\cdots
\gamma_{p}\beta_{p+1}^{-1}\gamma_{p}^{-1}\beta_{p}^{-1}\gamma_{p-1}^{-1}\cdots
\gamma_{1}^{-1}\beta_{1}^{-1}\alpha_{1}\]
\[=\gamma_{q-1}^{-1}\cdots
\gamma_{2}^{-1}\gamma_{q+2}\cdots
\gamma_{p}\gamma_{p}^{-1}\gamma_{p-1}^{-1}\cdots \gamma_{1}^{-1}\]
\[=\gamma_{q-1}^{-1}\cdots
\gamma_{2}^{-1}\gamma_{q+1}^{-1}\cdots \gamma_{1}^{-1}\]
\[=\alpha_{q}\alpha_{2}^{-1}\alpha_{q+2}
\alpha_{1}^{-1}\]
\[=\alpha_{q}\alpha_{q+2}=1\]\\ In general

\[b_{2m-1}^{2}=\beta_{1+q-m}^{-1}\gamma_{q-m}^{-1}\cdots
\gamma_{3}^{-1}\beta_{3}^{-1}\gamma_{2}^{-1}\beta_{2}\cdots
\beta_{q+m}\gamma_{q+m+1}\cdots
\gamma_{p}\beta_{p+1}^{-1}\gamma_{p}^{-1}\beta_{p}^{-1}\gamma_{p-1}^{-1}\cdots
\gamma_{1}^{-1}\beta_{1}^{-1}\alpha_{1}\]
\[=\gamma_{q-m}^{-1}\cdots
\gamma_{2}^{-1}\gamma_{q+m+1}\cdots
\gamma_{p}\gamma_{p}^{-1}\gamma_{p-1}^{-1}\cdots \gamma_{1}^{-1}\]
\[=\gamma_{q-m}^{-1}\cdots
\gamma_{2}^{-1}\gamma_{q+m}^{-1}\cdots \gamma_{1}^{-1}\]
\[=\alpha_{q-m+1}\alpha_{2}^{-1}\alpha_{q+m+1}\alpha_{1}^{-1}\]
\begin{equation} \label{alpha relations2}
=\alpha_{q-m+1}\alpha_{q+m+1}=1
\end{equation}
 $q-2\geq m\geq 1$\\

 The product relations in \ref{alpha relations1}
 and \ref{alpha relations2} can be summarized as

\begin{eqnarray*}
\alpha_{3}\alpha_{p+1}&=&1 \\
\alpha_{4}\alpha_{p}&=&1 \\
&\vdots& \\
\alpha_{q}\alpha_{q+3}&=&1 \\
\alpha_{q+1}\alpha_{q+2}&=&1
\end{eqnarray*}

and

\begin{eqnarray*}
\alpha_{q}\alpha_{q+2}&=&1 \\
\alpha_{q-1}\alpha_{q+3}&=&1 \\
&\vdots& \\
\alpha_{4}\alpha_{p-1}&=&1 \\
\alpha_{3}\alpha_{p}&=&1
\end{eqnarray*}

respectively. Using $c_{1}^{2}=\alpha _{q+1}=1$ we see that
$\alpha _{q+2}=1$. If $\alpha _{q+2}=1$, then  $\alpha _{q}=1$ and
if $\alpha _{q}=1$, then $\alpha _{q+3}=1$, etc. Continuing this
was we can conclude that $\alpha _{3}=\alpha _{4}=\cdots =\alpha
_{p+1}=1$ in the quotient $\pi_1\left(\Sigma_{p+1}\right)/N,$
which is isomorphic to $\pi_1\left(X\right)$,  and this concludes
the proof.

We do not know a topological argument to determine what those
simply connected symplectic $4-$ manifolds are but we can compute
their homeomorphism invariants using the symplectic Lefschetz
fibration structure that they support. This will be accomplished
by computing the Euler characteristics and the
signatures of those Lefschetz fibrations .\\

Another fact that we will borrow from the theory of Lefschetz
fibrations is the fact that the Euler characteristic of such a
fibration $X^{4}\rightarrow S^{2}$ is given by the formula
\[ \chi (X)=4-4g+s,\]\\ where $s$ is the number of singular
fibers, i.e., the number of vanishing cycles.\\

From \ref{theta p} we see that each of  $\theta_1^p$  and
$\theta_2^p$ consists of $p+9$ cycles. Therefore
$\phi_p^p=\left(\theta_2^p\theta_1^p\right)^p$ consists of
\[
2p\left( p+9\right)
\]\\ cycles and the Euler characteristic is found to be
\[ \chi (X)=4-4(p+1)+2p\left( p+9\right)
\]
\[ =2p^{2}+14p
\]
\[=2p\left( p+7\right).\]\\
In order to compute the signature $\sigma (X)$ we wrote a Matlab
program using the algorithm described in \cite{Oz}. The
computations that we have done using this program point to the
closed formula
\[\sigma (X)=-12p \]\\for the signature. Taking this formula for
granted we can conclude that
\[c_{1}^{2}\left( X\right) =3\sigma
\left( X\right) +2\chi \left( X\right) \] \[=3\left( -12p\right)
+2\cdot \allowbreak 2p\left( p+7\right) \] \[=4p\left( p-2\right)
\]\\ and
\[\chi _{h}\left( X\right) =\frac{\sigma \left( X\right) +\chi \left( X\right)
}{4}\]
\[=\frac{-12p+\allowbreak 2p\left( p+7\right) }{4}\]
\[=\frac{1}{2}p\left( p+1\right). \]\\

 $\chi_h(X)$ in the above computation makes sense because $X$ has
 almost complex structure.

\section{Appendix}

We include the actual outputs of signature computations for some
small values of $p$.\\

$p=3:$
\begin{eqnarray*}0 ,0 ,0 ,0 ,0 ,0 ,-1 ,-1 ,-1 ,-1
 ,-1 ,-1 ,-1 ,0 ,0 ,-1 ,-1 ,-1 ,-1 ,-1 ,0 ,0 ,0 ,0,\\
 0 ,0 ,0 ,0 ,-1 ,-1 ,-1 ,-1 ,-1 ,0 ,0 ,-1 ,-1 ,0 ,0 ,-1
 ,-1 ,-1 ,-1 ,-1 ,0 ,0 ,0 ,0,\\
0 ,0 ,0 ,0 ,-1 ,-1 ,-1 ,-1 ,-1 ,0 ,0 ,-1
 ,-1 ,-1 ,0 ,-1 ,-1 ,-1 ,-1 ,0 ,0 ,0 ,0 ,0
 \end{eqnarray*}\\

$p=5:$
\begin{eqnarray*}
0 ,0 ,0 ,0 ,0 ,0 ,0 ,0 ,-1 ,-1 ,-1 ,-1 ,-1 ,-1 ,-1 ,0 ,0 ,0 ,0 ,-1
-1 ,-1 ,-1 ,-1 ,0 ,0 ,0 ,0,\\
0 ,0 ,0 ,0 ,0 ,0 ,-1 ,-1 ,-1 ,-1 ,-1 ,0 ,0 ,-1 ,-1 ,0 ,0 ,0 ,0 ,-1
-1 ,-1 ,-1 ,-1 ,0 ,0 ,0 ,0,\\
0 ,0 ,0 ,0 ,0 ,0 ,-1 ,-1 ,-1 ,-1 ,-1 ,0 ,0 ,-1 ,-1 ,0 ,0 ,0 ,0 ,-1
-1 ,-1 ,-1 ,-1 ,0 ,0 ,0 ,0,\\
0 ,0 ,0 ,0 ,0 ,0 ,-1 ,-1 ,-1 ,-1 ,-1 ,0 ,0 ,-1 ,-1 ,0 ,0 ,0 ,0 ,-1
-1 ,-1 ,-1 ,-1 ,0 ,0 ,0 ,0,\\
0 ,0 ,0 ,0 ,0 ,0 ,-1 ,-1 ,-1 ,-1 ,-1 ,0 ,0 ,-1 ,-1 ,-1 ,0 ,0 ,0
,-1 ,-1 ,-1 ,-1 ,0 ,0 ,0 ,0 ,0
\end{eqnarray*}\\

$p=7:$
\begin{eqnarray*}
0 ,0 ,0 ,0 ,0 ,0 ,0 ,0 ,0 ,0 ,-1 ,-1 ,-1 ,-1 ,-1 ,-1 ,-1 ,0 ,0 ,0,
0 ,0 ,0 ,-1,-1 ,-1 ,-1 ,-1 ,0 ,0 ,0 ,0,\\
0 ,0 ,0 ,0 ,0 ,0 ,0 ,0 ,-1 ,-1 ,-1 ,-1 ,-1 ,0 ,0 ,-1 ,-1 ,0 ,0 ,0,
0 ,0 ,0 ,-1,-1 ,-1 ,-1 ,-1 ,0 ,0 ,0 ,0,\\
0 ,0 ,0 ,0 ,0 ,0 ,0 ,0 ,-1 ,-1 ,-1 ,-1 ,-1 ,0 ,0 ,-1 ,-1 ,0 ,0 ,0,
0 ,0 ,0 ,-1 ,-1 ,-1 ,-1 ,-1 ,0 ,0 ,0 ,0,\\
0 ,0 ,0 ,0 ,0 ,0 ,0 ,0 ,-1 ,-1 ,-1 ,-1 ,-1 ,0 ,0 ,-1 ,-1 ,0 ,0 ,0,
0 ,0 ,0 ,-1 ,-1 ,-1 ,-1 ,-1 ,0 ,0 ,0 ,0,\\
0 ,0 ,0 ,0 ,0 ,0 ,0 ,0 ,-1 ,-1 ,-1 ,-1 ,-1 ,0 ,0 ,-1 ,-1 ,0 ,0 ,0,
0 ,0 ,0 ,-1 ,-1 ,-1 ,-1 ,-1 ,0 ,0 ,0 ,0,\\
0 ,0 ,0 ,0 ,0 ,0 ,0 ,0 ,-1 ,-1 ,-1 ,-1 ,-1 ,0 ,0 ,-1 ,-1 ,0 ,0 ,0,
0 ,0 ,0 ,-1 ,-1 ,-1 ,-1 ,-1 ,0 ,0 ,0 ,0,\\
0 ,0 ,0 ,0 ,0 ,0 ,0 ,0 ,-1 ,-1 ,-1 ,-1 ,-1 ,0 ,0 ,-1 ,-1 ,-1 ,0
,0, 0 ,0 ,0 ,-1 ,-1 ,-1 ,-1 ,0 ,0 ,0 ,0 ,0
\end{eqnarray*}\\

$p=9:$
\begin{eqnarray*}
0 ,0 ,0 ,0 ,0 ,0 ,0 ,0 ,0 ,0 ,0 ,0 ,-1 ,-1 ,-1 ,-1 ,-1 ,-1 ,-1 ,0,
0 ,0 ,0 ,0 ,0 ,0 ,0 ,-1 ,-1 ,-1 ,-1 ,-1 ,0 ,0 ,0 ,0,\\
 0 ,0 ,0 ,0 ,0 ,0 ,0 ,0 ,0 ,0 ,-1 ,-1 ,-1 ,-1 ,-1 ,0 ,0 ,-1
 ,-1 ,0 ,0 ,0 ,0 ,0 ,0 ,0 ,0 ,-1 ,-1 ,-1 ,-1 ,-1 ,0 ,0 ,0 ,0,\\
 0 ,0 ,0 ,0 ,0 ,0 ,0 ,0 ,0 ,0 ,-1 ,-1 ,-1 ,-1 ,-1 ,0 ,0 ,-1
 ,-1 ,0 ,0 ,0 ,0 ,0 ,0 ,-1 ,0 ,-1 ,-1 ,-1 ,-1 ,0 ,0 ,0 ,0 ,0,\\
 0 ,0 ,0 ,0 ,0 ,0 ,0 ,0 ,0 ,0 ,-1 ,-1 ,-1 ,-1 ,-1 ,0 ,0 ,-1
 ,-1 ,0 ,0 ,0 ,0 ,0 ,0 ,0 ,0 ,-1 ,-1 ,-1 ,-1 ,-1 ,0 ,0 ,0 ,0,\\
 0 ,0 ,0 ,0 ,0 ,0 ,0 ,0 ,0 ,0 ,-1 ,-1 ,-1 ,-1 ,-1 ,0 ,0 ,-1
 ,-1 ,0 ,0 ,0 ,0 ,0 ,0 ,0 ,0 ,-1 ,-1 ,-1 ,-1 ,-1 ,0 ,0 ,0 ,0,\\
 0 ,0 ,0 ,0 ,0 ,0 ,0 ,0 ,0 ,0 ,-1 ,-1 ,-1 ,-1 ,-1 ,0 ,0 ,-1
 ,-1 ,0 ,0 ,0 ,0 ,0 ,0 ,-1 ,0 ,-1 ,-1 ,-1 ,-1 ,0 ,0 ,0 ,0 ,0,\\
  0 ,0 ,0 ,0 ,0 ,0 ,0 ,0 ,0 ,0 ,-1 ,-1 ,-1 ,-1 ,-1 ,0 ,0 ,-1
 ,-1 ,0 ,0 ,0 ,0 ,0 ,0 ,0 ,0 ,-1 ,-1 ,-1 ,-1 ,-1 ,0 ,0 ,0 ,0,\\
 0 ,0 ,0 ,0 ,0 ,0 ,0 ,0 ,0 ,0 ,-1 ,-1 ,-1 ,-1 ,-1 ,0 ,0 ,-1
 ,-1 ,0 ,0 ,0 ,0 ,0 ,0 ,0 ,0 ,-1 ,-1 ,-1 ,-1 ,-1 ,0 ,0 ,0 ,0,\\
 0 ,0 ,0 ,0 ,0 ,0 ,0 ,0 ,0 ,0 ,-1 ,-1 ,-1 ,-1 ,-1 ,0 ,0 ,-1
 ,-1 ,-1 ,0 ,0 ,0 ,0 ,0 ,0 ,0 ,-1 ,-1 ,-1 ,-1 ,0 ,0 ,0 ,0 ,0
\end{eqnarray*}

\bibliographystyle{amsplain}

\end{document}